\documentclass[article,12pt]{amsart}

\usepackage{amsfonts}
\usepackage{amsmath}
\usepackage{mathrsfs}
\usepackage{graphicx}
\usepackage{color}
\usepackage{epsfig}
\usepackage{enumerate}

\numberwithin{equation}{section}
\theoremstyle{plain}
\newtheorem{thm}{Theorem}[section]
\newtheorem{rem}{Remark}[section]

\newtheorem{defi}{Definition}[section]

\newcommand{\dE}{\mathbb{E}}
\newcommand{\dR}{\mathbb{R}}
\newcommand{\dL}{\mathbb{L}}

\newcommand{\dP}{\mathbb{P}}

\newcommand{\dZ}{\mathbb{Z}}

\newcommand{\cN}{\mathcal{N}}

\newcommand{\rI}{\mathrm{I}}
\newcommand{\cF}{\mathcal{F}}

\newcommand{\veps}{\varepsilon}

\newcommand{\ind}{\mbox{1}\kern-.25em \mbox{I}}
\font\calcal=cmsy10 scaled\magstep1
\def\build#1_#2^#3{\mathrel{\mathop{\kern 0pt#1}\limits_{#2}^{#3}}}
\def\liml{\build{\longrightarrow}_{}^{{\mbox{\calcal L}}}}
\def\limp{\build{\longrightarrow}_{}^{\dP}}

\def\videbox{\mathbin{\vbox{\hrule\hbox{\vrule height1.4ex \kern.6em\vrule height1.4ex}\hrule}}}
\def\demend{\hfill $\videbox$\\}

\keywords{Elephant random walk, Multi-dimensional martingales, Almost sure convergence, Asymptotic normality}

\begin{document}
\title[On the multi-dimensional elephant random walk]
{On the multi-dimensional elephant random walk \vspace{1ex}}
\author{Bernard Bercu}
\address{Universit\'e de Bordeaux, Institut de Math\'ematiques de Bordeaux,
UMR 5251, 351 Cours de la Lib\'eration, 33405 Talence cedex, France.}
\author{Lucile Laulin}
\address{Ecole normale sup\'erieure de Rennes,
D\'epartement de Math\'ematiques,
Campus de Ker lann, Avenue Robert Schuman, 35170 Bruz
}
\thanks{}

\begin{abstract}
The purpose of this paper is to investigate the asymptotic behavior of the
multi-dimensional elephant random walk (MERW). It is a non-Markovian random walk 
which has a complete memory of its entire history. A wide range of literature is available
on the one-dimensional ERW. Surprisingly, no references are available on the MERW. The goal of this paper 
is to fill the gap by extending the results on the one-dimensional ERW to the MERW. 
In the diffusive and critical regimes, we establish the almost sure convergence, the law of iterated logarithm and 
the quadratic strong law for the MERW. The asymptotic normality of the MERW, properly normalized, is also provided. 
In the superdiffusive regime, we prove the almost sure convergence as well as the mean square convergence
of the MERW. All our analysis relies on asymptotic results for multi-dimensional martingales.
\end{abstract}
\maketitle

\ \vspace{-7ex}
\section{Introduction}
\label{S-I}


The elephant random walk (ERW) is a fascinating discrete-time random process arising from mathematical physics.
It is a non-Markovian random walk on $\dZ$ which has a complete memory of its entire history.
This anomalous random walk was introduced by Sch\"utz and Trimper \cite{Schutz04}, 
in order to investigate how long-range memory affects the random walk and induces a crossover
from a diffusive to superdiffusive behavior. It was referred to as the ERW in allusion to the traditional 
saying that elephants can always remember where they have been. The ERW shows three differents regimes 
depending on the location of its memory parameter $p$ which lies between zero and one.
\vspace{1ex}\\
Over the last decade, the ERW has received considerable attention in the mathematical physics literature
in the diffusive regime $p< 3/4$ and the critical regime $p=3/4$, see
e.g.\! \cite{Baur16},\cite{Boyer14},\cite{Cressoni13},\cite{Cressoni07},\cite{Da13},\cite{Kumar10},\cite{Kursten16},\cite{Par06}
and the references therein. 
Quite recently, Baur and Bertoin \cite{Baur16} and independently Coletti, Gava and Sch\"utz \cite{Col17}
have proven the asymptotic normality of the ERW, properly normalized, with an explicit asymptotic variance. 
\vspace{1ex}\\
The superdiffusive regime $p>3/4$ is much harder to handle.
Initially, it was suggested by Sch\"utz and Trimper \cite{Schutz04} that, even in the superdiffusive regime,
the ERW has a Gaussian limiting distribution. However, it turns out \cite{Bercu17} that this limiting 
distribution is not Gaussian, as it was already predicted in \cite{Da13}, see also 
\cite{Col17},\cite{Par06}.
\vspace{1ex}\\
Surprisingly, to the best of our knowledge, no references are available on the multi-dimensional elephant random walk
(MERW) on $\dZ^d$, except \cite{Cressoni13},\cite{Lyu17} in the special case $d=2$. 
The goal of this paper is to fill the gap by extending the results 
on the one-dimensional ERW to the MERW. To be more precise, we shall study the influence of the memory parameter $p$ on the
MERW and we will show that the critical value is given by
$$
p_d=\frac{2d+1}{4d}.
$$
In the diffusive and critical regimes $p\leq p_d$, the reader will find the natural extension to higher dimension 
of the results recently established in \cite{Baur16},\cite{Bercu17},\cite{Col17},\cite{ColN17} on the almost sure asymptotic behavior 
of the ERW as well as on its asymptotic normality. One can notice that unlike in the classic random walk, the asymptotic normality 
of the MERW holds in any dimension $d \geq 1$. In the superdiffusive regime $p>p_d$, we will also prove
some extensions of the results in \cite{Bercu17},\cite{Cressoni13},\cite{Lyu17}.
\vspace{1ex}\\
Our strategy is to make an extensive use of the theory of martinagles \cite{Duflo97},\cite{HallHeyde80}, in particular
the strong law of large numbers and the central limit theorem for multi-dimensional martingales \cite{Duflo97}, 
as well as the law of iterated logarithm \cite{Stout70},\cite{Stout74}.
We strongly believe that our approach could be successfully extended to MERW with stops \cite{Cressoni13},\cite{Harbola14}, 
to amnesiac MERW \cite{Cressoni07}, as well as to MERW with reinforced memory \cite{Baur16},\cite{Harris15}.
\vspace{1ex}\\
The paper is organized as follows. In Section \ref{S-MERW}, we introduce the exact MERW and 
the multi-dimensional martingale we will extensively make use of. The main results of the paper are given 
in Section \ref{S-MR}. As usual, we first investigate the diffusive regime $p<p_d$ and we establish the almost sure convergence, 
the law of iterated logarithm and the quadratic strong law for the MERW. The asymptotic normality of the MERW, 
properly normalized, is also provided. Next, we prove similar results in the critical regime $p=p_d$. 
At last, we study the superdiffusive regime $p>p_d$ and we prove the almost sure convergence as well as 
the mean square convergence of the MERW to a non-degenerate random vector. Our martingale approach is described 
in Appendix A, while all technical proofs are postponed to Appendices B and C.

\ \vspace{-2ex} \\
\section{The multi-dimensional elephant random walk}
\label{S-MERW}

First of all, let us introduce the MERW. It is the natural extension to higher dimension of the one-dimensional ERW
defined in the pioneer work of Sch\"utz and Trimper \cite{Schutz04}. For a given dimension $d \geq 1$, let $(S_n)$ be
a random walk on $\dZ^d$, starting at the origin at time zero, $S_0 = 0$. At time $n = 1$, the elephant moves 
in one of the $2d$ directions with the same probability $1/2d$.  Afterwards, at time $n \geq 1$, 
the elephant chooses uniformly at random an integer $k$ among the previous times $1,\ldots,n$. 
Then, he moves exactly in the same direction as that of time $k$ with probability $p$ or in one of the
$2d-1$ remaining directions with the same probability $(1-p)/(2d-1)$, where the parameter $p$ 
stands for the memory parameter of the MERW.
From a mathematical point of view, the step of the elephant at time $n \geq 1$ in given by 
\begin{equation}
\label{INCMERW}
X_{n+1}=A_n X_k
\end{equation}
where 
\begin{equation*}
   A_{n} = \left \{ \begin{array}{ccc}
    +I_d &\text{ with probability } & p \vspace{1ex}\\
    -I_d &\text{ with probability } & \frac{1-p}{2d-1} \vspace{1ex}\\
     +J_d &\text{ with probability } & \frac{1-p}{2d-1} \vspace{1ex}\\
    -J_d &\text{ with probability } & \frac{1-p}{2d-1} \vspace{1ex}\\
     & \vdots & \vspace{1ex}\\
      +J_d^{d-1} &\text{ with probability } & \frac{1-p}{2d-1} \vspace{1ex}\\
    -J_d^{d-1}  &\text{ with probability } & \frac{1-p}{2d-1} 
   \end{array} \nonumber \right.
   \vspace{2ex}
\end{equation*}
with
\begin{equation*}
I_d =
\begin{pmatrix}
   1 & 0 & \cdots & 0 \\
   0 & 1 & \cdots & 0 \\
   \vdots & \ddots & \ddots &\vdots \\
   0 & \cdots& 0 & 1 
\end{pmatrix}
\hspace{1cm}
\text{and}
\hspace{1cm}
J_d = 
\begin{pmatrix}
   0 & 1 & \cdots & 0 \\
   0 & 0 & \cdots & 0 \\
   \vdots & \ddots & \ddots &\vdots \\
   1 & \cdots& 0 & 0 
\end{pmatrix}.
\end{equation*}
One can observe that the permutation matrix $J_d$ satisfies $J_d^d=I_d$.
Therefore, the position of the elephant at time $n \geq 1$ is given by
\begin{equation} 
\label{POSMERW}
S_{n+1} = S_n + X_{n+1}. 
\end{equation}
It follows from our very definition of the MERW that 
at any $n \geq 1$, $X_{n+1} = A_n X_{b_n}$ where $A_n$ is the random
matrix described before while $b_n$ is a random variable uniformly distributed on $\{1, ... , n\}$. Moreover, as
$A_n$ and $b_n$ are conditionally independent, we clearly have 
$\dE\left[X_{n+1}|\cF_n \right] = \dE\left[A_n\right] \dE\left[X_{b_n}|\cF_n\right]$
where $\cF_n$ stands for the $\sigma$-algebra, $\cF_n=\sigma(X_1, \ldots,X_n)$. Hence, we can deduce
from the law of total probability that at any time $n \geq 1$, 
\begin{equation}
\label{CEX}
\dE\left[X_{n+1}|\cF_n \right] =  \frac{1}{n}\Bigl(\frac{2dp-1}{2d-1}\Bigr) S_n = \frac{a}{n} S_n
\hspace{1cm}\text{a.s.}
\end{equation}
where $a$ is the fundamental parameter of the MERW,
\begin{equation}
\label{DEFA}
a=\frac{2dp-1}{2d-1}.
\end{equation}
Consequently, we immediately obtain from \eqref{POSMERW} and \eqref{CEX} that for any $n \geq 1$,
\begin{equation}
\label{CES1}
\dE\left[S_{n+1}|\cF_n\right] = \gamma_n S_n \hspace{1cm} \text{where} \hspace{1cm}
\gamma_n = 1+\frac{a}{n}.
\end{equation}
Furthermore, 
\begin{equation*}\prod_{k=1}^{n} \gamma_k=\frac{\Gamma(a+1+n)}{\Gamma(a+1)\Gamma(n+1)} 
\end{equation*}
where $\Gamma$ is the standard Euler Gamma function.
The critical value associated with the memory parameter $p$ of the MERW is
\begin{equation}
\label{DEFPD}
p_d=\frac{2d+1}{4d}.
\end{equation}
As a matter of fact,
$$
a<\frac{1}{2} \Longleftrightarrow p< p_d, 
\hspace{1cm}
a=\frac{1}{2} \Longleftrightarrow p= p_d,
\hspace{1cm}
a>\frac{1}{2} \Longleftrightarrow p> p_d.
$$
\begin{defi}
\label{DEF3R}
The MERW $(S_n)$ is said to be diffusive if $p<p_d$, critical if $p=p_d$, and superdiffusive if $p>p_d$.
\end{defi}
All our investigation in the three regimes relies on a martingale approach. To be more precise,
the asymptotic behavior of $(S_n)$ is closely related to the one of the
sequence $(M_n)$ defined, for all $n \geq 0$, by $M_n = a_nS_n$ where $a_0=1$,
$a_1 = 1$ and, for all $n \geq 2$,
\begin{equation}
\label{DEFAN}
a_n=\prod_{k=1}^{n-1}\gamma_k^{-1} = \frac{\Gamma(a+1)\Gamma(n)}{\Gamma(n+a)}. 
\end{equation}
It follows from a well-known property of the Euler Gamma function that
\begin{equation}
\label{CVGAMMAN}
\lim_{n \rightarrow \infty} \frac{\Gamma(n+a)}{\Gamma(n) n^a}= 1.
\end{equation}
Hence, we obtain from \eqref{DEFAN} and \eqref{CVGAMMAN} that
\begin{equation}
\label{CVGAN}
\lim_{n \rightarrow \infty} n^a a_n= \Gamma(a+1).
\end{equation}
Furthermore, since $a_n = \gamma_na_{n+1}$, we can deduce from \eqref{CES1} that for all $n \geq 1$, 
$$\dE\left[M_{n+1}|\cF_n\right] = M_n \hspace{1cm}  \text{a.s.}$$
It means that $(M_n)$ is a multi-dimensional martingale. Our goal is to extend 
the results recently established in \cite{Bercu17} to MERW. 
One can observe that our approach is much more tricky than that of \cite{Bercu17} as it requires to study
the asymptotic behavior of the  multi-dimensional martingale $(M_n)$.

\section{Main results}
\label{S-MR}
\subsection{The diffusive regime} 
Our first result deals with the strong law of large numbers for the MERW in the diffusive regime where $0 \leq p < p_d$.

\begin{thm}
\label{T-ASCVG-DR}
We have the almost sure convergence
\begin{equation}
\label{T-ASCVG-DR1}
\lim_{n \to \infty} \frac{1}{n}S_n = 0 \hspace{1cm}\text{a.s.}
\end{equation}
\end{thm}

\noindent
Some refinements on the almost sure rates of convergence for the MERW are as follows.

\begin{thm}
\label{T-ASCVGRATES-DR}
We have the quadratic strong law
\begin{equation}
\label{T-ASCVG-DR2}
\lim_{n \to \infty} \frac{1}{\log n}\sum_{k=1}^n \frac{1}{k^2}S_kS_k^T=\frac{1}{d(1-2a)}I_d \hspace{1cm} \text{a.s.}
\end{equation}
In particular,
\begin{equation}
\label{T-ASCVG-DR3}
\lim_{n \to \infty} \frac{1}{\log n}\sum_{k=1}^n \frac{\|S_k\|^2}{k^2}=\frac{1}{(1-2a)} \hspace{1cm} \text{a.s.}
\end{equation}
Moreover, we also have the law of iterated logarithm
\begin{equation}
\label{T-ASCVG-DR4}
 \limsup_{n \rightarrow \infty} \frac{\|S_n\|^2}{2 n \log \log n}
  = \frac{1}{(1-2a)} \hspace{1cm} \text{a.s.}
\end{equation}
\end{thm}

\noindent
Our next result is devoted to the asymptotic normality of the MERW in the diffusive regime $0 \leq p <p_d$.

\begin{thm}
\label{T-AN-DR}
We have the asymptotic normality
\begin{equation}
\label{T-CLT-DR}
\frac{1}{\sqrt{n}} S_n \liml \cN \Bigl(0, \frac{1}{(1-2a)d} I_d\Bigr).
\end{equation}
\end{thm}

\begin{rem} 
We clearly have from \eqref{DEFA} that
$$
\frac{1}{1-2a}=\frac{2d-1}{2d(1-2p)+1}.
$$
Hence, in the special case $d=1$, the critical value $p_d=3/4$ and the asymptotic variance
$$
\frac{1}{1-2a}=\frac{1}{3-4p}.
$$
Consequently, we find again the asymptotic normality for the one-dimensional ERW 
in the diffusive regime $0 \leq p <3/4$ recently established in \cite{Baur16},\cite{Bercu17},\cite{Col17}.
\end{rem}

\subsection{The critical regime}
We now focus our attention on the critical regime where the memory parameter $p = p_d$.
\begin{thm}
\label{T-ASCVG-CR}
We have the almost sure convergence
\begin{equation}
\label{T-ASCVG-CR1}
\lim_{n \to \infty} \frac{1}{\sqrt{n}\log n}S_n = 0 \hspace{1cm} \text{a.s.}
\end{equation}
\end{thm}

\noindent
We continue with some refinements on the almost sure rates of convergence for the MERW.

\begin{thm}
\label{T-ASCVGRATES-CR}
We have the quadratic strong law
\begin{equation}
\label{T-ASCVG-CR2}
\lim_{n \rightarrow \infty}  \frac{1}{\log \log n}\sum_{k=2}^n\frac{1}{(k \log k)^2}S_kS_k^T=\frac{1}{d}I_d
\hspace{1cm}  \text{a.s}
\end{equation}
In particular,
\begin{equation}
\label{T-ASCVG-CR3}
\lim_{n \to \infty}  \frac{1}{\log \log n}\sum_{k=2}^n\frac{\|S_k\|^2}{(k \log k)^2}=1 \hspace{1cm} \text{a.s.}
\end{equation}
Moreover, we also have the law of iterated logarithm
\begin{equation}
\label{T-ASCVG-CR4}
\limsup_{n \rightarrow \infty} \frac{\|S_n\|^2}{2 n\log n\log \log \log n}=1  \hspace{1cm}  \text{a.s.}
\end{equation}
\end{thm}

\noindent
Our next result concerns the asymptotic normality of the MERW in the critical regime $p=p_d$.

\begin{thm}
\label{T-AN-CR}
We have the asymptotic normality
\begin{equation}
\label{T-AN-CR1}
\frac{1}{\sqrt{n\log n}} S_n\liml \cN \Bigl(0, \frac{1}{d} I_d\Bigr).
\end{equation}
\end{thm}

\begin{rem} 
As before, in the special case $d=1$, we find again \cite{Baur16},\cite{Bercu17},\cite{Col17}
the asymptotic normality for the one-dimensional ERW 
$$
\frac{S_n}{\sqrt{n\log n}} \liml \cN (0, 1).
$$
\end{rem}

\subsection{The superdiffusive regime}
Finally, we get a handle on the more arduous superdiffusive regime where $p_d < p \leq 1$.

\begin{thm}
\label{T-ASCVG-SR}
We have the almost sure convergence
\begin{equation}
\label{T-ASCVG-SR1}
\lim_{n \to \infty} \frac{1}{n^a} S_n = L  \hspace{1cm}  \text{a.s.}
\end{equation}
where the limiting value L is a non-degenerate random vector. Moreover, we also have the mean square convergence
\begin{equation}
\label{T-ASCVG-SR2}
\lim_{n \to \infty} \dE\Bigl[\Bigl\|\frac{1}{n^a}S_n-L\Bigr\|^2\Bigr]=0.
\end{equation}
\end{thm}

\begin{thm}
\label{T-MOM-SR}
The expected value of  $L$ is $\dE[L]=0$, while its covariance matrix
is given by
\begin{equation}
\label{COVL}
\dE\left[LL^T\right]=\frac{1}{d(2a -1)\Gamma(2a)} I_d.
\end{equation}
In particular,
\begin{equation}
\label{TRCOVL}
\dE\left[\|L\|^2\right]=\frac{1}{(2a-1)\Gamma(2a)}.
\end{equation}
\end{thm}

\begin{rem}
Another possibility for the MERW is that, at time $n = 1$, the elephant moves in one direction, 
say the first direction $e_1$ of the standard basis $(e_1, \ldots, e_d)$ of $\dR^d$, with probability
$q$ or in one of the $2d-1$ remaining directions with the same probability $(1-q)/(2d-1)$, where 
the parameter $q$ lies in the interval $[0,1]$. Afterwards, at any time $n \geq 2$, the elephant
moves exactly as before, which means that his steps are given by \eqref{INCMERW}.
Then, the results of Section \ref{S-MR} holds true except Theorem \ref{T-MOM-SR} where
$$
\dE[L]=\frac{1}{\Gamma(a+1)}\Bigl( \frac{2dq -1}{2d-1} \Bigr) e_1
$$
and
$$
\dE[LL^T]= \frac{1}{\Gamma(2a+1)} \Bigl( \frac{2dq -1}{2d-1} \Bigr)\Bigl(e_1e_1^T - \frac{1}{d}I_d\Bigr)
+ \frac{1}{d(2a -1)\Gamma(2a)} I_d,
$$
which also leads to
$$
\dE\left[\|L\|^2\right]=\frac{1}{(2a-1)\Gamma(2a)}.
$$
\end{rem}

\section*{Appendix A \\ A multi-dimensional martingale approach}
\renewcommand{\thesection}{\Alph{section}}
\renewcommand{\theequation}{\thesection.\arabic{equation}}
\setcounter{section}{1}
\setcounter{equation}{0}

We clearly obtain from \eqref{INCMERW} that for any time $n\geq 1$, $\|X_{n} \| = 1$. Consequently,
it follows from \eqref{POSMERW} that $\|S_n\| \leq n$. Therefore, the sequence $(M_n)$ given, for all $n \geq 0$, by $M_n = a_nS_n$,
is a locally square-integrable multi-dimensional martingale. It can be rewritten in the additive form
\begin{equation}
\label{DEFMN}
M_n = \sum_{k=1}^{n}a_k \veps_k
\end{equation}
since its increments $\Delta M_n = M_n - M_{n-1}$ satisfy $\Delta M_n = a_nS_n - a_{n-1}S_{n-1} = a_n \veps_n$
where $\veps_n = S_n - \gamma_{n-1}S_{n-1}$. The predictable quadratic variation associated with
$(M_n)$ is the random square matrix of order $d$ given, for all $n \geq 1$, by
\begin{equation}
\label{IPMN}
\langle M \rangle_n = \sum_{k=1}^{n} \dE \left[\Delta M_k (\Delta M_k)^T | \cF_{k-1}\right].
\end{equation}
We already saw from \eqref{CES1} that $\dE\left[\veps_{n+1}|\cF_n\right] = 0$. Moreover, we deduce from
\eqref{POSMERW} together with \eqref{CEX} that
\begin{eqnarray}
\dE\left[S_{n+1} S_{n+1}^T|\cF_n\right] 
& = & \dE\left[S_nS_n^T|\cF_n\right]+\frac{2a}{n}S_nS_n^T + \dE\left[X_{n+1}X_{n+1}^T|\cF_n\right] \notag\\
& = & \left(1+ \frac{2a}{n}\right)S_nS_n^T + \dE\left[X_{n+1}X_{n+1}^T|\cF_n\right]\hspace{1cm} \text{a.s.}
\label{CES2}						
\end{eqnarray}
In order to calculate the right-hand side of \eqref{CES2}, one can notice that for any $n \geq 1$,
$$X_{n}X_{n}^T = \sum_{i=1}^d \rI_{X^i_{n} \neq 0}e_ie_i^{T} $$ 
where $(e_1, \ldots, e_d)$ stands for the standard basis of the Euclidean space $\dR^d$ and $X^i_{n}$ is the
$i$-th coordinate of the random vector $X_n$.
Moreover, it follows from \eqref{INCMERW} together with the law of total probability that any time $n \geq 1$
and for any $1\leq i \leq d$,
\begin{eqnarray*}
\dP(X_{n+1}^i\neq 0 | \cF_n)& = &\frac{1}{n}\sum _{k=1}^n \dP((A_nX_{k})^i\neq0 | \cF_n)  \\
& = &\frac{1}{n}\sum _{k=1}^n \rI_{X_k^i \neq 0} \dP(A_n = \pm I_d)+
\frac{1}{n}\sum _{k=1}^n
(1-\rI_{X_k^i \neq 0}) \dP(A_n = \pm J_d)  \\
& = &\frac{N_n^{X}(i)}{n} \Bigl(\dP(A_n=I_d)-\dP(A_n=J_d)\Bigr)+ 2\dP(A_n=J_d) 
\end{eqnarray*}
which implies that for any $1\leq i \leq d$,
\begin{equation}
\label{EXNNEQ}
\dE\bigl[\rI_{X_{n+1}^i \neq 0}|\cF_n\bigr] = \frac{a}{n} N_n^{X}(i) +
\frac{(1-a)}{d}  \hspace{1cm} \text{a.s.}
\end{equation}
where 
$$
N_n^{X}(i) = \sum_{k=1}^n \rI_{X_k^i \neq 0}
$$
and the parameter $a$ is given by \eqref{DEFA}.
Hence, we infer from \eqref{CES2} and \eqref{EXNNEQ} that
\begin{equation}
\label{CEX2}
\dE\left[X_{n+1}X_{n+1}^T|\cF_n\right] = \frac{a}{n}\Sigma_n + \frac{(1-a)}{d}I_d \hspace{1cm} \text{a.s.}
\end{equation}
where 
\begin{equation}
\label{DEFSIGMAN}
\Sigma_n = \sum_{i=1}^d N_n^{X}(i) e_ie_i^{T}. 
\end{equation}
One can observe the elementary fact that for all $n\geq 1$, $\text{Tr}(\Sigma_n)=n$
where $\text{Tr}(\Sigma_n)$ stands for the trace of the positive definite matrix $\Sigma_n$.
Therefore, we obtain from \eqref{CES2} together with \eqref{CEX2} that 
\begin{eqnarray}
\dE\left[\veps_{n+1}\veps_{n+1}^T|\cF_n\right] & = & \dE\left[S_{n+1}S_{n+1}^T|\cF_n\right]-\gamma_n^2S_nS_n^T \notag \\
& = & \Bigl(1+ \frac{2a}{n}\Bigr)S_nS_n^T + \frac{a}{n}\Sigma_n + \frac{(1-a)}{d}I_d -\gamma_n^2S_nS_n^T \notag \\
& = & \frac{a}{n}\Sigma_n + \frac{(1-a)}{d}I_d - \left( \frac{a}{n}\right) ^2 S_nS_n^T\hspace{1cm}  \text{a.s.}
\label{CEMEPS2}
\end{eqnarray} 
which ensures that
\begin{eqnarray}
\label{CEEPS2}
\dE\left[\|\veps_{n+1}\|^2|\cF_n\right]
& = &\frac{a}{n}\text{Tr}(\Sigma_n) + \frac{1-a}{d}\text{Tr}(I_d) - \left( \frac{a}{n}\right) ^2 \|S_n\|^2 \notag \\
& = & 1 - ( \gamma_n -1 ) ^2 \|S_n\|^2 \hspace{1cm}\text{a.s.}
\end{eqnarray}				
By the same token, 
$$
\dE\left[\|\veps_{n+1}\|^4|\cF_n\right] = 1 - 3(\gamma_n -1)^4\|S_n\|^4  -2 (\gamma_n-1)^2\|S_n\|^2 + 4(\gamma_n-1)^2\xi_n
$$
where, thanks to \eqref{CEX2},
$$
\xi_n  = \dE\left[\langle S_n , X_{n+1}\rangle ^2 | \cF_n\right] =\frac{a}{n}S_n^T\Sigma_nS_n + \frac{(1-a)}{d}\|S_n\|^2.
$$
It leads to
\begin{eqnarray}
\dE\left[\|\veps_{n+1}\|^4|\cF_n\right] 
& =  &  1 - 3(\gamma_n -1)^4\|S_n\|^4  - 2 \Bigl(1-\frac{2(1-a)}{d}\Bigr)(\gamma_n-1)^2\|S_n\|^2 \notag \\
&  & \hspace{5ex} +\frac{4a}{n}(\gamma_n-1)^2S_n^T\Sigma_nS_n \hspace{1cm} \text{a.s.}
\label{CEEPS4}
\end{eqnarray}
Therefore, as $\Sigma_n \leq nI_d$ for the usual order of positive definite matrices, we clearly obtain from
\eqref{CEEPS4} that 
\begin{eqnarray}
\label{MAJCEEPS4}
\dE\left[\|\veps_{n+1}\|^4|\cF_n\right] & \leq & 1 - 3(\gamma_n -1)^4\|S_n\|^4  \notag \\
& &\hspace{4ex} +\frac{2}{d} (\gamma_n-1)^2
\Bigl(2a(d-1) +2-d\Bigr)\|S_n\|^2\hspace{0.5cm} \text{a.s.}
\end{eqnarray}
Consequently, we obtain from \eqref{CEEPS2} and \eqref{MAJCEEPS4} the almost sure upper bounds
\begin{equation}
\label{MAJMOMEPS}
\sup_{n\geq0} \dE\left[\|\veps_{n+1}\|^2|\cF_n\right] \leq 1 \hspace{1cm} 
\text{and} \hspace{1cm} 
\sup_{n\geq0} \dE\left[\|\veps_{n+1}\|^4|\cF_n\right] \leq \frac{4}{3} \hspace{1cm} \text{a.s.}
\end{equation}			
Hereafter, we deduce from \eqref{IPMN} and \eqref{CEMEPS2} that
\begin{eqnarray}
\langle M \rangle_n 
& = & a_1^2 \dE[\veps_1\veps_1^T] + \sum_{k=1}^{n-1} a_{k+1}^2 \dE\left[\veps_{k+1}\veps_{k+1}^T|\cF_k\right] \notag \\
& = & \frac{1}{d}I_d \sum_{k=1}^{n}a_k^2 + a \sum_{k=1}^{n-1}a_{k+1}^2 \Bigl( \frac{1}{k}\Sigma_k-\frac{1}{d}I_d\Bigr) - \zeta_n
\label{CALIPMN}
\end{eqnarray}
where
$$
\zeta_n=a^2 \sum_{k=1}^{n-1}\Bigl(\frac{a_{k+1}}{k}\Bigr)^2 S_kS_k^T.
$$
Hence, by taking the trace on both sides of \eqref{CALIPMN}, we find that
\begin{equation}
\label{TRIPMN}
\text{Tr}\langle M \rangle_n =\sum_{k=1}^{n}a_k^2 - a ^2 \sum_{k=1}^{n-1} \Bigl(\frac{a_{k+1}}{k}\Bigr)^2 \| S_k \|^2.
\end{equation}
The asymptotic behavior of the multi-dimensional martingale $(M_n)$ is closely related to the one of
\begin{equation*}
v_n=\sum_{k=1}^na_k^2=\sum_{k=1}^n\Bigl(\frac{\Gamma(a+1)\Gamma(k)}{\Gamma(a+k)}\Bigr)^2.
\end{equation*} 
One can observe that we always have $\text{Tr} \langle M \rangle_n \leq v_n$.
In accordance with Definition \ref{DEF3R}, we have three regimes. In the diffusive regime where $a<1/2$,
\begin{equation}
\label{VNDIFF}
\lim_{n \to \infty} \frac{v_n}{n^{1-2a}}= \ell \hspace{1cm} \text{where} \hspace{1cm}  
\ell=\frac{(\Gamma(a+1))^2}{1-2a}.
\end{equation}
In the critical regime where $a=1/2$,
\begin{equation}
\label{VNCRIT}
\lim_{n\to \infty} \frac{v_n}{\log n} = (\Gamma(a+1))^2=\frac{\pi}{4}.
\end{equation}
Finally, in the superdiffusive regime where $a>1/2$, $v_n$ converges to the finite value
\begin{eqnarray}
\label{VNSUPER}
\lim_{n\rightarrow \infty} v_n 
& = &\sum_{k=0}^{\infty}  \Bigl( \frac{\Gamma(a+1) \Gamma(k+1)  }{\Gamma(a+k+1)} \Bigr)^2
= \sum_{k=0}^{\infty} \frac{(1)_k\,(1)_k\,(1)_k} {(a+1)_k\, (a+1)_k\,k!} \notag\\
& = & {}_{3}F_2\Bigl( \begin{matrix}
{1,1,1}\\
{a+1,a+1}\end{matrix} \Bigl| 1\Bigr)
\end{eqnarray}
where, for any $\alpha\in \dR$, $(\alpha)_k=\alpha(\alpha+1)\cdots(\alpha+k-1)$ for $k\geq 1$, 
$(\alpha)_0=1$ stands for the Pochhammer symbol and ${}_{3}F_2$ is the
generalized hypergeometric function defined by
\begin{eqnarray*}
{}_{3}F_2 \Bigl( \begin{matrix}
{a,b,c}\\
{d,e}\end{matrix} \Bigl|
{\displaystyle z}\Bigr)
=\sum_{k=0}^{\infty}
\frac{(a)_k\,(b)_k\,(c)_k}
{(d)_k\,(e)_k\, k!} z^k.
\end{eqnarray*}

\section*{Appendix B \\ Proofs of the almost sure convergence results}
\renewcommand{\thesection}{\Alph{section}}
\renewcommand{\theequation}{\thesection.\arabic{equation}}
\setcounter{section}{2}
\setcounter{equation}{0}

\vspace{-3ex}
\subsection*{}
\begin{center}
{\bf B.1. The diffusive regime.}
\end{center}
\ \vspace{-4ex}\\

\noindent{\bf Proof of Theorem \ref{T-ASCVG-DR}.}
First of all, we focus our attention on the proof of the
almost sure convergence \eqref{T-ASCVG-DR1}. We already saw from 
\eqref{TRIPMN} that $\text{Tr} \langle M \rangle_n \leq v_n$. Moreover, we obtain from \eqref{VNDIFF} that, in the diffusive regime
where $0<a<1/2$, $v_n$ increases to infinity with the speed $n^{1-2a}$. On the one hand, it follows from the strong law of large numbers for multi-dimensional martingales given e.g. by the last part of Theorem 4.3.15 in \cite{Duflo97} that for any $\gamma>0$,
\begin{equation}
\label{SLLNMN}
\frac{\|M_n\|^2}{\lambda_{max} \langle M \rangle_n }=o\Bigl( \Bigl(\log \text{Tr} \langle M \rangle_n \Bigr)^{1+\gamma}\Bigr)
\hspace{1cm} \text{a.s}
\end{equation}
where $\lambda_{max} \langle M \rangle_n $ stands for the maximal eigenvalue of
the random square matrix $\langle M \rangle_n$. However, as $\langle M \rangle_n$ is a 
positive definite matrix and $\text{Tr} \langle M \rangle_n  \leq v_n$,
we clearly have $\lambda_{max} \langle M \rangle_n \leq \text{Tr} \langle M \rangle_n \leq v_n$.
Consequenly, we obtain from \eqref{SLLNMN} that
\begin{equation*}
\|M_n\|^2=o\bigl( v_n (\log v_n )^{1+\gamma}\bigr)
\hspace{1cm} \text{a.s}
\end{equation*}
which implies that
\begin{equation}
\label{CVGMNDIFF}
\|M_n\|^2 =o\bigl( n^{1-2a}(\log n)^{1+\gamma} \bigr) \hspace{1cm} \text{a.s.}
\end{equation}
Hence, as $M_n = a_nS_n$, it follows from \eqref{CVGAN} and \eqref{CVGMNDIFF} that
for any $\gamma>0$,
\begin{equation*}
\|S_n\|^2 =o\bigl(n (\log n)^{1+\gamma}\bigr) \hspace{1cm}  \text{a.s.}
\end{equation*}
which completes the proof of Theorem \ref{T-ASCVG-DR}.
\demend

\noindent{\bf Proof of Theorem \ref{T-ASCVGRATES-DR}.}
We shall now proceed to the proof of the almost sure rates of convergence given in 
Theorem \ref{T-ASCVGRATES-DR}. First of all, we claim that
\begin{equation}
\label{CVGSIGMAN}
\lim_{n\to\infty}\frac{1}{n}\Sigma_n = \frac{1}{d}I_d \hspace{1cm} \text {a.s.}
\end{equation}
where $\Sigma_n$ is the random square matrix of order $d$ given by \eqref{DEFSIGMAN}. As a matter of fact,
in order to prove \eqref{CVGSIGMAN}
it is only necessary to show that for any $1 \leq i \leq d$,
\begin{equation}
\label{CVGNXNI}
\lim_{n\to\infty}\frac{N_n^{X}(i)}{n} = \frac{1}{d} \hspace{1cm} \text {a.s.}
\end{equation}
For any $1 \leq i \leq d$, denote
$$\Lambda_n(i) = \frac{N_n^{X}(i)}{n}.$$
One can observe that 
$$\Lambda_{n+1}(i)=\frac{n}{n+1}\Lambda_n(i)+ \frac{1}{n+1}\rI_{X^i_{n+1} \neq 0}$$ 
which leads, via \eqref{EXNNEQ}, to the recurrence relation  
\begin{equation}
\label{RECLNI}
\Lambda_{n+1}(i) = \frac{n}{n+1}\gamma_n \Lambda_n(i) +  \frac{(1-a)}{d(n+1)} +\frac{1}{n+1}\delta_{n+1}(i)
\end{equation} 
where $\delta_{n+1}(i)=\rI_{X^i_{n+1} \neq 0}-\dE[\rI_{X^i_{n+1}\neq 0}|\cF_n]$.
After straightforward calculations, the solution of this recurrence relation is given by
\begin{equation}
\label{SOLLNI}
\Lambda_n(i) =\frac{1}{na_n}\Bigl(\Lambda_1(i) + \frac{(1-a)}{d} \sum_{k=2}^na_k + L_n(i) \Bigr)
\end{equation}
where
$$ 
L_n(i) = \sum_{k=2}^na_k \delta_k(i).
$$
However, $(L_n(i))$ is a square-integrable real martingale with predictable quadratic variation
$\langle L(i) \rangle_n$ satisfying $\langle L(i) \rangle_n \leq v_n$ a.s. Then, it follows from
the standard strong law of large numbers for martingales given by 
Theorem 1.3.24 in \cite{Duflo97} that $(L_n(i))^2=O(v_n \log v_n)$ a.s. Consequently, 
as $na_n^2$ is equivalent to $(1-2a) v_n$, we obtain that
for any $1 \leq i \leq d$,
\begin{equation}
\label{CVGLNI}
\lim_{n\to\infty}\frac{1}{na_n}L_n(i) = 0 \hspace{1cm} \text {a.s.}
\end{equation}
Furthermore, one can easily check from \eqref{CVGAN} that
\begin{equation}
\label{CVGSUMAN}
\lim_{n\to\infty}\frac{1}{na_n}\sum_{k=1}^na_k = \frac{1}{1-a}.
\end{equation}
Therefore, we find from \eqref{SOLLNI} together with \eqref{CVGLNI} and \eqref{CVGSUMAN} that
for any $1 \leq i \leq d$,
\begin{equation}
\label{CVGLAMBDANI}
\lim_{n\to\infty}\Lambda_n(i) = \frac{1}{d} \hspace{1cm} \text {a.s.}
\end{equation}
which immediately leads to \eqref{CVGNXNI}. Hereafter, it follows from the conjunction
of \eqref{T-ASCVG-DR1}, \eqref{CEMEPS2} and \eqref{CVGNXNI} that
\begin{equation}
\label{CVGEPS2}
\lim_{n\to\infty}
\dE\left[\veps_{n+1}\veps_{n+1}^T|\cF_n\right]= \frac{1}{d}I_d \hspace{1cm} \text {a.s.}
\end{equation}
By the same token, we also obtain from \eqref{CALIPMN} and Toeplitz lemma that
\begin{equation}
\label{CVGIPMN}
\lim_{n\to\infty}
\frac{1}{v_n} \langle M \rangle_n= \frac{1}{d}I_d \hspace{1cm} \text {a.s.}
\end{equation}
We are now in the position to prove the quadratic strong law 
\eqref{T-ASCVG-DR2}. For any vector $u$ of $\dR^d$, denote $M_n(u)=\langle u, M_n \rangle$
and $\veps_n(u)=\langle u, \veps_n \rangle$.
We clearly have from \eqref{DEFMN} 
\begin{equation*}
M_n(u) = \sum_{k=1}^{n}a_k \veps_k(u).
\end{equation*}
Consequently, $(M_n(u))$ is a square-integrable real martingale. Moreover, it follows from
\eqref{CVGEPS2} that
\begin{equation*}
\lim_{n\to\infty}
\dE\left[|\veps_{n+1}(u)|^2|\cF_n\right]= \frac{1}{d}\|u\|^2 \hspace{1cm} \text {a.s.}
\end{equation*}
Moreover, we can deduce from \eqref{MAJMOMEPS} and the Cauchy-Schwarz inequality that
\begin{equation*}
\sup_{n\geq 0} \dE\left[|\veps_{n+1}(u)|^4|\cF_n\right] \leq \frac{4}{3}\|u\|^4 \hspace{1cm} \text{a.s.}
\end{equation*}	
Furthermore, we clearly have from \eqref{CVGAN} and \eqref{VNDIFF} that
$$
\lim_{n \rightarrow \infty} n f_n= 1-2a \hspace{1cm} \text{where} \hspace{1cm} f_n=\frac{a_n^2}{v_n},
$$
which of course implies that $f_n$ converges to zero.
Therefore, it follows from the quadratic strong law for real martingales given e.g. in Theorem 3 of \cite{Bercu04},
that for any vector $u$ of $\dR^d$,
\begin{equation}
\label{LFQ-DR1}
\lim_{n\rightarrow \infty} 
\frac{1}{\log v_n} \sum_{k=1}^{n} f_{k} \Bigl( \frac{M_{k}^2(u)}{v_k} \Bigr)
= \frac{1}{d}\|u\|^2 \hspace{1cm} \text{a.s.}
\end{equation}
Consequently, we find from \eqref{VNDIFF} and \eqref{LFQ-DR1} that
\begin{equation}
\label{LFQ-DR2}
\lim_{n\rightarrow \infty} 
\frac{1}{\log n} \sum_{k=1}^{n}   \frac{a_k^2}{v_k^2} M_{k}^2(u) = \frac{(1-2a)}{d}\|u\|^2 \hspace{1cm} \text{a.s.}
\end{equation}
Hereafter, as $M_n=a_nS_n$ and $n^2 a_n^4$ is equivalent to $(1-2a)^2v_n^2$, we obtain
from \eqref{LFQ-DR2} that for any vector $u$ of $\dR^d$,
\begin{equation}
\label{LFQ-DR3}
\lim_{n\rightarrow \infty} 
\frac{1}{\log n} \sum_{k=1}^n  \frac{1}{k^2}u^TS_kS_k^T u=\frac{1}{d(1-2a)}\|u\|^2 \hspace{1cm} \text{a.s.}
\end{equation}
By virtue of the second part of Proposition 4.2.8 in \cite{Duflo97},
we can conclude from \eqref{LFQ-DR3} that
\begin{equation}
\label{LFQ-DR4}
\lim_{n\rightarrow \infty} 
\frac{1}{\log n} \sum_{k=1}^n  \frac{1}{k^2}S_kS_k^T =\frac{1}{d(1-2a)}I_d \hspace{1cm} \text{a.s.}
\end{equation}
which completes the proof of \eqref{T-ASCVG-DR2}. By taking the trace on both sides of \eqref{LFQ-DR4},
we immediately obtain \eqref{T-ASCVG-DR3}.
Finally, we shall proceed to the proof of the law of iterated logarithm
given by \eqref{T-ASCVG-DR4}. We already saw that $a_n^4v_n^{-2}$ is equivalent to $(1-2a)^2n^{-2}$.
It ensures that
\begin{equation}
\label{CONDLIL}
\sum_{n=1}^{+\infty} \frac{a_{n}^4}{v_n^2} < +\infty. 
\end{equation}
Hence, it follows from the law of iterated logarithm for real martingales due to
Stout \cite{Stout70},\cite{Stout74}, see also Corollary 6.4.25 in \cite{Duflo97}, that for any vector $u$ of $\dR^d$,
\begin{eqnarray}
 \limsup_{n \rightarrow \infty} \Bigl(\frac{1}{2 v_n \log \log v_n}\Bigr)^{1/2} M_n(u) & = & 
 -\liminf_{n \rightarrow \infty} \Bigl(\frac{1}{2 v_n \log \log v_n}\Bigr)^{1/2} M_n(u) \nonumber \\
 & = & \frac{1}{\sqrt{d}}\|u\| \hspace{1cm} \text{a.s.}
 \label{LIL-MG-DR}
\end{eqnarray}
Consequently, as $M_n(u)=a_n\langle u, S_n\rangle$, we obtain from \eqref{VNDIFF} together with \eqref{LIL-MG-DR} that
\begin{eqnarray*}
 \limsup_{n \rightarrow \infty} \Bigl(\frac{1}{2 n \log \log n}\Bigr)^{1/2} \langle u, S_n\rangle & = & 
 -\liminf_{n \rightarrow \infty} \Bigl(\frac{1}{2 n \log \log n}\Bigr)^{1/2} \langle u, S_n\rangle \nonumber \\
 & = & \frac{1}{\sqrt{d(1-2a)}}\|u\| \hspace{1cm} \text{a.s.}
\end{eqnarray*}
In particular, for any vector $u$ of $\dR^d$,
\begin{equation}
\label{LIL-SNU-DR}
 \limsup_{n \rightarrow \infty} \frac{1}{2 n \log \log n}\langle u, S_n\rangle^2
  = \frac{1}{d(1-2a)}\|u\|^2 \hspace{1cm} \text{a.s.}
\end{equation}
However,
$$
\| S_n \|^2=\sum_{i=1}^d \langle e_i, S_n\rangle^2
$$
where $(e_1, \ldots, e_d)$ is the standard basis of $\dR^d$. Finally, we deduce from
\eqref{LIL-SNU-DR} that
\begin{equation*}
 \limsup_{n \rightarrow \infty} \frac{\|S_n\|^2}{2 n \log \log n}
  = \frac{1}{(1-2a)} \hspace{1cm} \text{a.s.}
\end{equation*}
which achieves the proof of Theorem \ref{T-ASCVGRATES-DR}.
\demend

\vspace{-2ex}

\subsection*{}
\begin{center}
{\bf B.2. The critical regime.}
\end{center}
\ \vspace{-4ex}\\

\noindent{\bf Proof of Theorem \ref{T-ASCVG-CR}.}
We already saw from \eqref{VNCRIT} that in the critical regime where $a=1/2$, $v_n$ increases slowly to
infinity with a logarithmic speed $\log n$. We obtain once again from the last part of Theorem 4.3.15 in \cite{Duflo97}
that for any $\gamma>0$,
\begin{equation*}
\|M_n\|^2=o\bigl( v_n (\log v_n )^{1+\gamma}\bigr)
\hspace{1cm} \text{a.s}
\end{equation*}
which leads to
\begin{equation}
\label{CVGMNCRIT}
\|M_n\|^2 =o\bigl( \log n (\log\log n)^{1+\gamma} \bigr) \hspace{1cm} \text{a.s.}
\end{equation}
However, we clearly have from \eqref{CVGAN} with $a=1/2$ that
\begin{equation}
\label{CVGANCRIT}
\lim_{n\to \infty}na_n^2=\frac{\pi}{4}.
\end{equation}
Consequently, as $M_n = a_nS_n$, we deduce from \eqref{CVGMNCRIT} and \eqref{CVGANCRIT} that
for any $\gamma>0$,
\begin{equation*}
\|S_n\|^2 =o\bigl(n \log n (\log \log n)^{1+\gamma}\bigr) \hspace{1cm}  \text{a.s.}
\end{equation*}
which completes the proof of Theorem \ref{T-ASCVG-CR}.
\demend

\vspace{-2ex}


\noindent{\bf Proof of Theorem \ref{T-ASCVGRATES-CR}.} The proof of Theorem
\ref{T-ASCVGRATES-CR} is left to the reader as it follows the same lines as that of
Theorem \ref{T-ASCVGRATES-DR}.
\demend
\vspace{-2ex}

\subsection*{}
\begin{center}
{\bf B.3. The superdiffusive regime.}
\end{center}
\ \vspace{-4ex} \\

\noindent{\bf Proof of Theorem \ref{T-ASCVG-SR}.}
We already saw from \eqref{VNSUPER} that in the superdiffusive regime where $1/2<a \leq 1$, 
$v_n$ converges to a finite value. 
As previously seen, $\text{Tr} \langle M \rangle_n \leq v_n$. Hence, we clearly have
$$
\lim_{n \rightarrow \infty} \text{Tr} \langle M \rangle_n < \infty \hspace{1cm} \text{a.s.}
$$
Therefore, if
\begin{equation}
\label{DEFLNSUPER}
L_n=\frac{M_n}{\Gamma(a+1)},
\end{equation}
we can deduce from the second part of Theorem 4.3.15 in \cite{Duflo97} that
\begin{equation}
\label{CVGMNSUPER}
\lim_{n \to \infty} M_n = M \hspace{1cm} \text{and} \hspace{1cm} \lim_{n \to \infty} L_n = L \hspace{1cm} \text{a.s.}
\end{equation}
where the limiting values $M$ and $L$ are the random vectors of $\dR^d$ given by
$$
M=\sum_{k=1}^{\infty}a_k\veps_k \hspace{1cm} \text{and} \hspace{1cm} L=\frac{1}{\Gamma(a+1)}\sum_{k=1}^{\infty}a_k\veps_k.
$$
Consequently, as $M_n = a_nS_n$, \eqref{T-ASCVG-SR1} clearly follows from \eqref{CVGAN} and \eqref{CVGMNSUPER}
We now focus our attention on the mean square convergence \eqref{T-ASCVG-SR2}. As $M_0=0$, we have
from \eqref{DEFMN} and \eqref{IPMN} that for all $n \geq 1$,
$$
\dE[\|M_n\|^2]= \sum_{k=1}^n \dE[\|\Delta M_{k}\|^2]=\dE[\text{Tr} \langle M \rangle_n] \leq v_n.
$$
Hence, we obtain from \eqref{VNSUPER} that
\begin{equation*}
\sup_{n \geq 1} \dE\left[\|M_n\|^2\right] \leq 
{}_{3}F_2\Bigl( \begin{matrix}
{1,1,1}\\
{a+1,a+1}\end{matrix} \Bigl| 1\Bigr)< \infty,
\end{equation*}
which means that the martingale $(M_n)$ is bounded in $\dL^2$. 
Therefore, we have the mean square convergence
\begin{equation*}
\lim_{n \to \infty} \dE\bigl[\|M_n-M\|^2\bigr]=0,
\end{equation*}
which clearly leads to \eqref{T-ASCVG-SR2}.
\demend

\noindent{\bf Proof of Theorem \ref{T-MOM-SR}.}
First of all, we clearly have for all $n \geq 1$, $\dE[M_n]=0$ which implies that $\dE[M]=0$ leading to $\dE[L]=0$.
Moreover, taking expectation on both sides of \eqref{CES2} and \eqref{CEX2}, we obtain 
that for all $n \geq 1$,
\begin{eqnarray}
\dE\left[S_{n+1}S_{n+1}^T\right] & = & 
\Bigl( 1+\frac{2a}{n}\Bigr)\dE\left[S_{n}S_n^T\right] + \dE \left[X_{n+1}X_{n+1}^T\right]
\notag	\\ 
& = & \Bigl( 1+\frac{2a}{n}\Bigr)\dE\left[S_{n}S_n^T\right] +\frac{a}{n}\dE\left[\Sigma_n\right]+\frac{(1-a)}{d}I_d.
\label{COVSN-1}
\end{eqnarray}
However, we claim that
\begin{equation}
\label{ESIGMAN}
\dE\left[\Sigma_n\right] = \frac{n}{d}I_d.
\end{equation}
As a matter of fact, taking expectation on both sides of \eqref{SOLLNI}, we find that
for any $1 \leq i \leq d$,
\begin{equation}
\label{ELNI}
\dE[\Lambda_n(i)] =\frac{1}{na_n}\Bigl(\dE[\Lambda_1(i)] + \frac{(1-a)}{d} \sum_{k=2}^na_k \Bigr).
\end{equation}
On the one hand, we clearly have 
$$\dE[\Lambda_1(i)]=\frac{1}{d}.$$
On the other hand, it follows from Lemma B.1 in \cite{Bercu17} that
\begin{eqnarray}
\label{SUMAN}
\sum_{k=2}^na_k  & = & \sum_{k=2}^n \frac{\Gamma(a+1)\Gamma(k)}{\Gamma(k+a)}
= \sum_{k=1}^{n-1} \frac{\Gamma(a+1)\Gamma(k+1)}{\Gamma(k+a+1)} \notag \\
& = & \frac{1}{(a-1)}
\left(1 - \frac{\Gamma(a+1)\Gamma(n+1)}{\Gamma(a+n)}\right) = 
\frac{(1 - na_n)}{(a-1)}.
\end{eqnarray}
Consequently, we can deduce from \eqref{ELNI} and \eqref{SUMAN} that for any $1 \leq i \leq d$,
\begin{equation}
\label{ELNIFIN}
\dE[\Lambda_n(i)] =\frac{1}{na_n}\Bigl(\frac{1}{d} - \frac{(1-na_n)}{d}  \Bigr)=\frac{1}{d}.
\end{equation}
Therefore, we get from \eqref{DEFSIGMAN} and \eqref{ELNIFIN} that
\begin{equation*}
\dE[\Sigma_n] = n \sum_{i=1}^d \dE[\Lambda_n(i)] e_ie_i^{T}= \frac{n}{d} 
\sum_{i=1}^d  e_ie_i^{T}=\frac{n}{d} I_d.
\end{equation*}
Hereafter, we obtain from \eqref{COVSN-1} and \eqref{ESIGMAN} that
\begin{equation}
\dE\left[S_{n+1}S_{n+1}^T\right] =\Bigl( 1+\frac{2a}{n}\Bigr)\dE\left[S_{n}S_n^T\right]+ 
\frac{1}{d} I_d.
\label{COVSN-2}
\end{equation}
It is not hard to see that the solution of this recurrence relation is given by
\begin{eqnarray}
\label{COVSN-3}
\dE\left[S_{n}S_n^T\right]  & = & \frac{\Gamma(n+2a)}{\Gamma(2a+1)\Gamma(n)}
\left( \dE[S_1S_1^T] + \frac{1}{d} \sum_{k=1}^{n-1} \frac{\Gamma(2a+1)\Gamma(k+1)}{\Gamma(k+2a+1)} I_d \right) \notag \\
& = & \frac{\Gamma(n+2a)}{\Gamma(n)}
\left(  \sum_{k=1}^{n} \frac{\Gamma(k)}{\Gamma(k+2a)}  \right) \frac{1}{d}I_d
\end{eqnarray}
since
$$
\dE[S_1S_1^T] = \frac{1}{d}I_d.
$$
Therefore, it follows once again from Lemma B.1 in \cite{Bercu17} that
\begin{equation}
\label{COVSN-4}
\dE\left[S_{n}S_n^T\right]=\frac{n}{(2a -1)}\left( \frac{\Gamma(n+2a)}{\Gamma(n+1)\Gamma(2a)} -1\right) \frac{1}{d}I_d.
\end{equation}
Hence, we obtain from \eqref{DEFLNSUPER} together with \eqref{COVSN-4} that
\begin{eqnarray}
\label{COVLN}
\dE[L_n L_n^T] & = & \frac{na_n^2}{(2a -1)(\Gamma(a+1))^2}\left( \frac{\Gamma(n+2a)}{\Gamma(n+1)\Gamma(2a)} -1\right) \frac{1}{d}I_d \notag \\
& = & \frac{n}{(2a -1)} \left( \frac{\Gamma(n)}{\Gamma(n+a)}\right)^2 \left( \frac{\Gamma(n+2a)}{\Gamma(n+1)\Gamma(2a)} -1\right) \frac{1}{d}I_d.
\end{eqnarray}
Finally, we find from \eqref{T-ASCVG-SR2} and \eqref{COVLN} that
$$
\lim_{n\rightarrow \infty} \dE[L_n L_n^T]=\dE[L L^T]=
\frac{1}{d(2a -1)\Gamma(2a)} I_d
$$
which achieves the proof of Theorem \ref{T-MOM-SR}.
\demend

\section*{Appendix C \\ Proofs of the asymptotic normality results}
\renewcommand{\thesection}{\Alph{section}}
\renewcommand{\theequation}{\thesection.\arabic{equation}}
\setcounter{section}{3}
\setcounter{equation}{0}

\vspace{-2ex}
\subsection*{}
\begin{center}
{\bf C.1. The diffusive regime.}
\end{center}
\ \vspace{-4ex}\\

\noindent{\bf Proof of Theorem \ref{T-AN-DR}.}
In order to establish the asymptotic normality \eqref{T-CLT-DR}, we shall make use of the central limit theorem for 
multi-dimensional martingales given e.g. by Corollary 2.1.10 of \cite{Duflo97}.
First of all, we already saw from \eqref{CVGIPMN} that
\begin{equation}
\label{RCVGIPMN}
\lim_{n\to\infty}
\frac{1}{v_n} \langle M \rangle_n= \frac{1}{d}I_d \hspace{1cm} \text {a.s.}
\end{equation}
Consequently, it only remains to show that $(M_n)$ satisfies Lindeberg's condition, in other words, for all $\veps > 0$,
\begin{equation*}
\frac{1}{v_n}\sum_{k=1}^n \dE\left[\| \Delta M_n \|^2 \rI_{\{\| \Delta M_n \| \geq \veps \sqrt{v_n}\}}| \cF_{k-1}\right] \limp 0.
\end{equation*}
We have from \eqref{MAJMOMEPS} that for all $\veps > 0$
\begin{eqnarray*}
\frac{1}{v_n}\sum_{k=1}^n \dE\left[\| \Delta M_n \|^2 \rI_{\{\| \Delta M_n \| \geq \veps \sqrt{v_n}\}}| \cF_{k-1}\right] 
& \leq &\frac{1}{\veps^2 v_n^2}\sum_{k=1}^n \dE\left[\| \Delta M_n \|^4 | \cF_{k-1}\right] \\
& \leq & \sup_{1 \leq k \leq n} \dE\left[\|\veps_k\|^4|\cF_{k-1}\right]\frac{1}{\veps^2 v_n^2} \sum_{k=1}^na_k^4 \\
& \leq & \frac{4}{3\veps^2 v_n^2}\sum_{k=1}^na_k^4.
\end{eqnarray*}
However, we already saw from \eqref{CONDLIL} that
\begin{equation*}
\sum_{n=1}^{+\infty} \frac{a_{n}^4}{v_n^2} < +\infty. 
\end{equation*}
Hence, it follows from Kronecker's lemma that 
$$\lim_{n\to \infty}\frac{1}{v_n^2}\sum_{k=1}^na_k^4 =0,$$
which ensures that Lindeberg's condition is satisfied. Therefore, we can conclude
from the central limit theorem for martingales that
\begin{equation}
\label{CLTMN-DR}
\frac{1}{\sqrt{v_n}}M_n\liml \mathcal{N}\Bigl(0,\frac{1}{d}I_d\Bigr).
\end{equation}
As $M_n = a_nS_n$ and $\sqrt{n}a_n$ is equivalent to $\sqrt{v_n(1-2a)}$, we find from \eqref{CLTMN-DR} that
$$\frac{1}{\sqrt{n}} S_n\liml \mathcal{N}\Bigl(0,\frac{1}{d(1-2a)}I_d\Bigr),$$
which completes the proof of Theorem \ref{T-AN-DR}.
\demend

\vspace{-5ex}

\subsection*{}
\begin{center}
{\bf C.2. The critical regime.}
\end{center}
\ \vspace{-4ex}\\

\noindent{\bf Proof of Theorem \ref{T-AN-CR}.}
Via the same lines as in the proof of \eqref{CVGIPMN}, we can deduce
from \eqref{T-ASCVG-CR1}, (\ref{TRIPMN}) and \eqref{VNCRIT} that in the critical regime
\begin{equation}
\lim_{n \to \infty} \frac{1}{v_n} \langle M\rangle_n=\frac{1}{d}I_d \hspace{1cm} \text{a.s.}
\end{equation}
Moreover, it follows from \eqref{VNCRIT} and \eqref{CVGANCRIT} that $a_n^2v_n^{-1}$ is equivalent to $(n\log n)^{-1}$.
It implies that
\begin{equation}
\label{CONDLILCR}
\sum_{k=1}^{\infty}\frac{a_n^4}{v_n^2} < + \infty.
\end{equation}
As previously seen, we infer from \eqref{CONDLILCR} that $(M_n)$
satisfies Lindeberg's condition. Therefore, we can conclude from the central limit
theorem for martingales that
\begin{equation}
\label{CLTMN-CR}
\frac{1}{\sqrt{v_n}} M_n\liml \mathcal{N}\Bigl(0,\frac{1}{d}I_d\Bigr).
\end{equation}
Finally, as $M_n = a_nS_n$ and $a_n\sqrt{n \log n}$ is equivalent to $\sqrt{v_n}$, we obtain from that
\eqref{CLTMN-CR} that
\begin{equation*}
\frac{1}{\sqrt{n \log n}} S_n \liml \cN (0,1),
\end{equation*}
which achieves the proof of Theorem \ref{T-AN-CR}.
\demend

\vspace{-2ex}

\bibliographystyle{acm}

\end{document}